\input amsppt.sty
\hsize=12.5cm
\vsize=540pt

\def\Aut{\operatorname{Aut}}
\def\CC{\Bbb C}
\def\CO{\Cal O}
\def\dist{\operatorname{dist}}
\def\eps{\varepsilon}
\def\id{\operatorname{id}}
\def\NN{\Bbb N}
\def\RR{\Bbb R}
\def\too{\longrightarrow}

\redefine\kappa{\varkappa}
\redefine\phi{\varphi}
\redefine\theta{\vartheta}

\NoRunningHeads

\topmatter
\title Kobayashi--Royden vs\. Hahn pseudometric in $\CC^2$\endtitle
\author Witold Jarnicki\endauthor

\address
\noindent
Instytut Matematyki\newline
Uniwersytet Jagiello\'nski\newline
Reymonta 4\newline
30-059 Krak\'ow, Poland
\endaddress

\email wmj\@im.uj.edu.pl \endemail

\abstract
For a domain $D\subset\CC$ the Kobayashi--Royden $\kappa$ and Hahn $h$ 
pseudometrics are equal iff $D$ is simply connected. Overholt showed that for 
$D\subset\CC^n$, $n\geq3$, we have $h_D\equiv\kappa_D$. Let 
$D_1,D_2\subset\CC$. The aim of this paper is to show that 
$h_{D_1\times D_2}\equiv\kappa_{D_1\times D_2}$ iff at least one of $D_1$, 
$D_2$ is simply connected or biholomorphic to $\CC\setminus\{0\}$. In 
particular, there are domains $D\subset\CC^2$ for which 
$h_D\not\equiv\kappa_D$.
\endabstract

\endtopmatter

\document

\subhead
1. Introduction
\endsubhead

For a domain $D\subset\CC^n$, the Kobayashi--Royden pseudometric
$\kappa_D$ and the Hahn pseudometric $h_D$ are defined by the formulas:
$$
\align
\kappa_D(z;X)&:=\inf\{|\alpha|\:\exists_{f\in\CO(E,D)}\;f(0)=z,\;\alpha f'(0)=X\},\\
h_D(z;X)&:=\inf\{|\alpha|\:\exists_{f\in\CO(E,D)}\;f(0)=z,\;\alpha f'(0)=X,\text{ $f$ is injective}\},\\
&\hskip240pt z\in D,X\in\CC^n,
\endalign
$$
where $E$ denotes the unit disc (cf\. \cite{Roy}, \cite{Hah},
\cite{Jar-Pfl}). Obviously $\kappa_D\leq h_D$. It is known that both
pseudometrics are invariant under biholomorphic mappings, i.e., if
$f\:D\too\widetilde D$ is biholomorphic, then
$$
\multline
h_D(z;X)=h_{\widetilde D}(f(z);f'(z)(X)),\quad
\kappa_D(z;X)=\kappa_{\widetilde D}(f(z);f'(z)(X)),\\
z\in D,X\in\CC^n.
\endmultline
$$
It is also known that for a domain $D\subset\CC$ we have:
$h_D\equiv\kappa_D$ iff $D$ is simply connected. In particular
$h_D\not\equiv\kappa_D$ for $D=\CC_\ast:=\CC\setminus\{0\}$. It has turned
out that $h_D\equiv\kappa_D$ for any domain $D\subset\CC^n$, $n\geq3$
(\cite{Ove}). The case
$n=2$ was investigated for instance in \cite{Hah}, \cite{Ves}, \cite{Vig},
\cite{Cho}, but neither a proof nor a counterexample for the equality was
found (existing `counterexamples' were based on incorrect product properties
of the Hahn pseudometric).

\subhead
2. The main result
\endsubhead

\proclaim{Theorem 1} Let $D_1,D_2\subset\CC$ be domains. Then:

{\rm 1.} If at least one of $D_1$, $D_2$ is simply connected, then
$h_{D_1\times D_2}\equiv\kappa_{D_1\times D_2}$.

{\rm 2.} If at least one of $D_1$, $D_2$ is biholomorphic to $\CC_\ast$, then
$h_{D_1\times D_2}\equiv\kappa_{D_1\times D_2}$.

{\rm 3.} Otherwise $h_{D_1\times D_2}\not\equiv\kappa_{D_1\times D_2}$.
\endproclaim

Let $p_j\:D_j^\ast\too D_j$ be a holomorphic universal covering of $D_j$
($D_j^\ast\in\{\CC,E\}$), $j=1,2$. Recall that if $D_j$ is simply connected,
then $h_{D_j}\equiv\kappa_{D_j}$. If $D_j$ is not simply connected and $D_j$
is not biholomorphic to $\CC_\ast$, then, by the uniformization theorem,
$D_j^\ast=E$ and $p_j$ is not injective.

Hence, Theorem 1 is an immediate consequence of the following three
propositions (we keep the above notation).

\proclaim{Proposition 2} If $h_{D_1}\equiv\kappa_{D_1}$, then
$h_{D_1\times D_2}\equiv\kappa_{D_1\times D_2}$ for any domain
$D_2\subset\CC$.
\endproclaim

\proclaim{Proposition 3} If $D_1$ is biholomorphic to $\CC_\ast$, then
$h_{D_1\times D_2}\equiv\kappa_{D_1\times D_2}$ for any domain
$D_2\subset\CC$.
\endproclaim

\proclaim{Proposition 4} If $D_j^\ast=E$ and $p_j$ is not injective, $j=1,2$,
then $h_{D_1\times D_2}\not\equiv\kappa_{D_1\times D_2}$.
\endproclaim

Observe the following property that will be helpful in proving the propositions.
\proclaim{Remark 5} For any domain $D\subset\CC^n$ we have $h_D\equiv\kappa_D$ 
iff for any $f\in\CO(E,D)$, $\theta\in(0,1)$ with $f'(0)\neq0$, there exists an 
injective $g\in\CO(E,D)$ such that $g(0)=f(0)$ and $g'(0)=\theta f'(0)$.
\endproclaim

\demo{Proof of Proposition 2}
Let $f=(f_1,f_2)\in\CO(E,D_1\times D_2)$ and let $\theta\in(0,1)$.

First, consider the case where $f_1'(0)\neq0$.

By Remark 5, there exists an injective function $g_1\in\CO(E,D_1)$ such that
$g_1(0)=f_1(0)$ and $g_1'(0)=\theta f_1'(0)$. Put
$g(z):=(g_1(z),f_2(\theta z))$.

Obviously $g\in\CO(E,D_1\times D_2)$ and $g$ is injective. Moreover,
$g(0)=f(0)$ and
$g'(0)=(g_1'(0),f_2'(0)\theta)=(\theta f_1'(0),\theta f_2'(0))=\theta f'(0)$.

\medskip

Suppose now that $f_1'(0)=0$. Take $0<d<\dist(f_1(0),\partial D_1)$ and put
$$
\gather
h(z):=\frac{f_2(\theta z)-f_2(0)}{f_2'(0)},\quad
M:=\max\{|h(z)|\:z\in\overline E\},\\
g_1(z):=f_1(0)+\frac d{M+1}(h(z)-\theta z),\quad
g(z):=(g_1(z),f_2(\theta z)),\quad z\in E.
\endgather
$$
Obviously $g\in\CO(E,\CC\times D_2)$. Since $|g_1(z)-f_1(0)|<d$, we get
$g_1(z)\in B(f_1(0),d)\subset D_1$, $z\in E$. Hence
$g\in\CO(E,D_1\times D_2)$. Take $z_1,z_2\in E$ such that $g(z_1)=g(z_2)$.
Then $h(z_1)=h(z_2)$, and consequently $z_1=z_2$.

Finally $g(0)=(g_1(0),f_2(0))=(f_1(0)+\frac d{M+1}h(0),f_2(0))=f(0)$ and
$g'(0)=(g_1'(0),\theta f_2'(0))=(\frac d{M+1}(h'(0)-\theta),\theta f_2'(0))=
\theta f'(0)$.
\hfill$\square$\enddemo

\demo{Proof of Proposition 3} We may assume that $D_1=\CC_\ast$ and
$D_2\neq\CC$.
Using Remark 5, let $f=(f_1,f_2)\in\CO(E,\CC_\ast\times D_2)$ and let
$\theta\in(0,1)$. Applying an appropriate automorphism of $\CC_\ast$, we may
assume that $f_1(0)=1$.

\medskip

For the case where $f_2'(0)=0$, we apply the above construction to the
domains $\widetilde D_1=f_2(0)+\dist(f_2(0),\partial D_2)E$,
$\widetilde D_2=\CC_\ast$ and mappings $\widetilde f_1\equiv f_2(0)$,
$\widetilde f_2=f_1$.

\medskip

Now, consider the case where $f_2'(0)\neq0$ and $\theta f_1'(0)=1$. We put
$$
g_1(z):=1+z,\quad
g(z):=(g_1(z),f_2(\theta z)),\quad z\in E.
$$
Obviously, $g\in\CO(E,\CC_\ast\times D_2)$ and $g$ is injective. We have
$g(0)=(1,f_2(0))=f(0)$ and
$g'(0)=(1,\theta f_2'(0))=\theta f'(0)$.

\medskip

In all other cases, let $M:=\max\{|f_2(z)|\:|z|\leq\theta\}$. Take a
$k\in\NN$ such that $|c_k|>M$, where
$$
c_k:=f_2(0)-k\frac{\theta f_2'(0)}{\theta f_1'(0)-1}.
$$
Put
$$
\gather
h(z):=\frac{f_2(\theta z)-c_k}{f_2(0)-c_k},\\
g_1(z):=(1+z)h^k(z),\quad
g_2(z):=f_2(\theta z),\quad
g(z):=(g_1(z),g_2(z)),\quad z\in E.
\endgather
$$

Obviously, $g\in\CO(E,\CC\times D_2)$. Since $h(z)\neq0$, we
have $g_1(z)\neq0$, $z\in E$. Hence $g\in\CO(E,\CC_\ast\times D_2)$. Take
$z_1,z_2\in E$ such that $g(z_1)=g(z_2)$. Then $h(z_1)=h(z_2)$, and
consequently $z_1=z_2$.

Finally $g(0)=(h^k(0),f_2(0))=f(0)$ and
$$
\multline
g'(0)=(g_1'(0),\theta f_2'(0))=
(h^k(0)+kh^{k-1}(0)h'(0),\theta f_2'(0))\\=
\Big(1+k\frac{\theta f_2'(0)}{f_2(0)-c_k},\theta f_2'(0)\Big)=
(1+\theta f_1'(0)-1,\theta f_2'(0))=\theta f'(0).
\endmultline
$$

\hfill$\square$\enddemo

\demo{Proof of Proposition 4} It suffices to show that there exist
$\varphi_1,\varphi_2\in\Aut(E)$ and a point $q=(q_1,q_2)\in E^2$,
$q_1\neq q_2$, such that
$p_j(\phi_j(q_1))=p_j(\phi_j(q_2))$, $j=1,2$, and
$\det[(p_j\circ\varphi_j)'(q_k)]_{j,k=1,2}\neq0$.

Indeed, put $\widetilde p_j:=p_j\circ\varphi_j$, $j=1,2$, and suppose that
$h_{D_1\times D_2}\equiv\kappa_{D_1\times D_2}$. Put
$a:=(\widetilde p_1(0),\widetilde p_2(0))$ and
$X:=(\widetilde p_1'(0),\widetilde p_2'(0))\in(\CC_\ast)^2$.
Take an arbitrary $f\in\CO(E,D_j)$ with $f(0)=a_j$.
Let $\widetilde f$ be the lifting of $f$ with respect to $\widetilde p_j$
such that $\widetilde f(0)=0$. Since $|\widetilde f'(0)|\leq 1$, we get
$|f'(0)|\leq|X_j|$.
Consequently $\kappa_{D_j}(a_j;X_j)=1$, $j=1,2$. In particular,
$\kappa_{D_1\times D_2}(a;X)=
\max\{\kappa_{D_1}(a_1;X_1),\;\kappa_{D_2}(a_2;X_2)\}=1$.

Let
$(0,1)\ni\alpha_n\nearrow1$.
Fix an $n\in\NN$. Since $\kappa_{D_1\times D_2}(a;X)=1$, there exists
$f_n\in\CO(E,D_1\times D_2)$ such that $f_n(0)=a$ and $f_n'(0)=\alpha_n X$.
By Remark 5, there exists an injective holomorphic mapping
$g_n=(g_{n,1},g_{n,2})\:E\too D_1\times D_2$ such that $g_n(0)=a$ and
$g_n'(0)=\alpha_n^2X$. Let $\widetilde{g}_{n,j}$ be the lifting with respect
to $\widetilde p_j$ of $g_{n,j}$ with $\widetilde g_{n,j}(0)=0$, $j=1,2$.

By the Montel theorem, we may assume that the sequence
$(\widetilde g_{n,j})_{n=1}^\infty$ is locally uniformly convergent,
$\widetilde g_{0,j}:=\lim_{n\too\infty}\widetilde g_{n,j}$. We have
$\widetilde g'_{0,j}(0)=1$, $\widetilde g_{0,j}\:E\too E$. By the Schwarz
lemma we have $\widetilde g_{0,j}=\id_E$, $j=1,2$.

Let $h_{0,j}(z_1,z_2):=\widetilde p_j(z_1)-\widetilde p_j(z_2)$, $(z_1,z_2)\in E^2$,
$$
V_j=V(h_{0,j})=\{(z_1,z_2)\in E^2\:h_{0,j}(z_1,z_2)=0\},\quad j=1,2.
$$
Since
$$
\det\left[\frac{\partial h_{0,j}}{\partial z_k}(q)\right]_{j,k=1,2}
=-\det\left[\widetilde p_j'(q_k)\right]_{j,k=1,2}\neq0,
$$
$V_1$ and $V_2$ intersect
transversally at $q$. Let $U\subset\subset\{(z_1,z_2)\in E^2\:z_1\neq z_2\}$
be a neighborhood of $q$ such that
$V_1\cap V_2\cap\overline U=\{q\}$. For
$n\in\NN$, $j=1,2$, define
$$
h_{n,j}(z_1,z_2):=g_{n,j}(z_1)-g_{n,j}(z_2),\quad(z_1,z_2)\in E^2.
$$
Observe that the sequence $(h_{n,j})_{n=1}^\infty$ converges
uniformly on $\overline U$ to $h_{0,j}$, $j=1,2$. In particular
(cf\. \cite{Two-Win}), we have
$V(h_{n,1})\cap V(h_{n,2})\cap\overline U=
\{z\in\overline U\:h_{n,1}(z)=h_{n,2}(z)=0\}\neq\varnothing$ for some
$n\in\NN$ --- contradiction.

\medskip

We move now to the construction of $\varphi_1,\varphi_2$ and $q$.
Let $\psi_j\in\Aut(E)$ be a non--identity lifting of $p_j$ with respect to
$p_j$
($p_j\circ\psi_j\equiv p_j$, $\psi_j\not\equiv\id$), $j=1,2$. Observe that
$\psi_j$
has no fixed points (a lifting is uniquely determined by its value
at one point), $j=1,2$.

To simplify notation, let
$$
h_a(z):=\frac{z-a}{1-\overline az},\quad a,z\in E.
$$

One can easily check that
$$
\sup_{z\in E}m(z,\psi_j(z))=1,\quad j=1,2,
$$
where $m(z,w):=|h_w(z)|=\big|\frac{z-w}{1-z\overline w}\big|$ is the M\"obius
distance. Hence there exist $\eps\in(0,1)$ and $z_1,z_2\in E$ with
$m(z_1,\psi_1(z_1))=m(z_2,\psi_2(z_2))=1-\eps$. Let $d\in(0,1)$,
$h_1,h_2\in\Aut(E)$ be such that $h_j(-d)=z_j$, $h_j(d)=\psi_j(z_j)$, $j=1,2$.

If $(p_j\circ h_j)'(-d)\neq\pm(p_j\circ h_j)'(d)$ for
some $j$ (we may assume that for $j=1$), then at least one of the
determinants
$$
\gather
\det\left[\matrix
(p_1\circ h_1)'(-d),&(p_1\circ h_1)'(d)\\
(p_2\circ h_2)'(-d),&(p_2\circ h_2)'(d)
\endmatrix\right],\\
\det\left[\matrix
(p_1\circ h_1\circ(-\id))'(-d),&(p_1\circ h_1\circ(-\id))'(d)\\
(p_2\circ h_2)'(-d),&(p_2\circ h_2)'(d)
\endmatrix\right],
\endgather
$$
is nonzero.

Otherwise, let $\widetilde\psi_j=h_j^{-1}\circ\psi_j\circ h_j$ and
$\widetilde p_j=p_j\circ h_j$, $j=1,2$. Observe that $\widetilde\psi_j(-d)=d$
and $(\widetilde\psi_j'(-d))^2=1$, $j=1,2$. Thus, each $\widetilde\psi_j$ is
either $-\id$ or $h_c$, where $c=-\frac{2d}{1+d^2}$. The case
$\widetilde\psi_j=-\id$ is impossible since $\widetilde\psi_j$
has no fixed points. By substituting $p_j$ by
$\widetilde p_j$ and $\psi_j$ by $\widetilde\psi_j$, $j=1,2$, the proof
reduces to the
case, where $\psi_1=\psi_2=h_c=:\psi$ for some $-1<c<0$.

\medskip

We claim that there exists a
point $a\in E$ such that if an automorphism $\phi=\phi_a\in\Aut(E)$ satisfies
$\phi(a)=\psi(a)$ and $\phi(\psi(a))=a$,
then $\phi'(a)\neq\pm\psi'(a)$. Suppose for a moment that such an $a$ has
been found. Notice that $\phi\circ\phi=\id$ and hence
$\phi'(\psi(a))=\frac1{\phi'(a)}$.
Put $\varphi_1:=\id$, $\varphi_2:=\phi$, $q:=(a,\psi(a))$. We
have
$$
\align
&\det\left[\matrix
(p_1\circ\varphi_1)'(a),&(p_1\circ\varphi_1)'(\psi(a))\\
(p_2\circ\varphi_2)'(a),&(p_2\circ\varphi_2)'(\psi(a))
\endmatrix\right]\\
=&\det\left[\matrix
p_1'(a),&p_1'(\psi(a))\\
p_2'(\phi(a))\phi'(a),&p_2'(\phi(\psi(a))\phi'(\psi(a))
\endmatrix\right]\\
=&\det\left[\matrix
(p_1\circ\psi)'(a),&p_1'(\psi(a))\\
p_2'(\psi(a))\phi'(a),&(p_2\circ\psi)'(a)\frac1{\phi'(a)}
\endmatrix\right]\\
=&\det\left[\matrix
p_1'(\psi(a))\psi'(a),&p_1'(\psi(a))\\
p_2'(\psi(a))\phi'(a),&p_2'(\psi(a))\psi'(a)\frac1{\phi'(a)}
\endmatrix\right]\\
=&p_1'(\psi(a))p_2'(\psi(a))\det\left[\matrix
\psi'(a),&1\\
\phi'(a),&\frac{\psi'(a)}{\phi'(a)}
\endmatrix\right]\neq0,
\endalign
$$
which finishes the construction.

It remains to find $a$. First observe that the equality $\phi_a'(a)=\psi'(a)$
is impossible. Otherwise $\phi_a=\psi$ and consequently $\psi\circ\psi=\id$;
contradiction. We only need to find an $a\in E$ such that
$\phi_a'(a)\neq-\psi'(a)$. One can easily check that
$$
\phi_a=h_{-a}\circ(-\id)\circ h_{h_a(\psi(a))}\circ h_a.
$$
Direct calculations show that
$\phi_a'(a)=-\psi'(a)\Longleftrightarrow a\in\RR$. Thus it suffices to take
any $a\in E\setminus\RR$.
\hfill$\square$\enddemo

\subhead
3. Acknowledgement
\endsubhead

I would like to thank Professors Peter Pflug and W\l odzimierz Zwonek for
their valuable remarks.

\Refs\nofrills{References}

{

\widestnumber{\key}{Two-Win}
\ref
\key Cho
\by K\. S\. Choi
\paper Injective hyperbolicity of product domain
\jour J\. Korea Soc\. Math\. Educ\. Ser\. B Pure Appl\. Math\.
\vol 5(1)
\yr 1998
\pages 73--78
\endref
\ref
\key Hah
\by K\. T\. Hahn
\paper Some remark on a new pseudo-differential metric
\jour Ann\. Polon\. Math\.
\vol 39
\yr 1981
\pages 71--81
\endref
\ref
\key Jar-Pfl
\by M\. Jarnicki, P\. Pflug
\book Invariant Distances and Metrics in Complex Analysis
\publ de Gruyter Exp\. Math\. 9, Walter de Gruyter, Berlin
\yr 1993
\endref
\ref
\key Ove
\by M\. Overholt
\paper Injective hyperbolicity of domains
\jour Ann\. Polon\. Math\.
\vol 62(1)
\yr 1995
\pages 79--82
\endref
\ref
\key Roy
\by H\. L\. Royden
\book Remarks on the Kobayashi metric {\rm in ``Several complex
variables, II''}
\publ Lecture Notes in Math\. 189, Springer Verlag
\yr 1971
\pages 125--137
\endref
\ref
\key Two-Win
\by P\. Tworzewski, T\. Winiarski
\paper Continuity of intersection of analytic sets
\jour Ann\. Polon\. Math\.
\vol 42
\yr 1983
\pages 387--393
\endref
\ref
\key Ves
\by E\. Vesentini
\paper Injective hyperbolicity
\jour Ricerche Mat\.
\vol 36
\yr 1987
\pages 99--109
\endref
\ref
\key Vig
\by J.-P\. Vigu\'e
\paper Une remarque sur l'hyperbolicit\'e injective
\jour Atti\. Acc\. Lincei Rend\. fis\. (8)
\vol 83
\yr 1989
\pages 57--61
\endref

}

\endRefs

\enddocument